\documentclass{elsarticle}

\usepackage[all]{xy}

\usepackage{amsmath,amsthm}
\usepackage{amsfonts}
\usepackage{verbatim}
\usepackage{amssymb}
\usepackage{amscd}
\usepackage{bbm}
\usepackage{mathrsfs}
\usepackage{graphicx}

\newcommand{\set}[1]{\,\left\{#1\right\}}
\newcommand{\setd}[2]{\,\left\{#1\ \colon\ #2\right\}}

\newtheorem{theorem}{Theorem}
\newtheorem{corollary}[theorem]{Corollary}

\newtheorem{proposition}[theorem]{Proposition}
\newtheorem{definition}[theorem]{Definition}

\newtheorem{remark}[theorem]{Remark}
\newtheorem{example}[theorem]{Example}

\newtheorem{problem}[theorem]{Problem}

\newcommand{\diam}{\delta}

\newcommand{\sign}{\operatorname{sign}}
\newcommand{\disto}{c}
\newcommand{\edisto}{c_{(2)}}

\newcommand{\RR}{\mathbb{R}}
\newcommand{\ZZ}{\mathbb{Z}}

\newcommand{\NN}{\mathbb{N}}

\begin{document}

\begin{frontmatter}
\title{Diameters, distortion and eigenvalues}

\author{Rostislav I. Grigorchuk\fnref{g1}}
\ead{grigorch@math.tamu.edu}
\address{Department of Mathematics, Texas A\&M University, College Station, TX 77843}

\author{Piotr W. Nowak\fnref{n1}}
\ead{pnowak@math.tamu.edu}
\address{Department of Mathematics, Texas A\&M University, College Station, TX 77843}

\fntext[g1]{Partially supported by  NSF grant DMS-0600975}
\fntext[n1]{Partially supported by NSF grant DMS-0900874.}
\address{\vskip10pt\normalfont Dedicated to Antonio Mach\`{i}}

\begin{abstract}
We study the relation between the diameter, the first positive eigenvalue of 
the discrete $p$-Laplacian and the $\ell_p$-distortion of a finite graph. We prove 
an inequality relating 
these three quantities and apply it to families of Cayley and Schreier graphs. We also
show that the $\ell_p$-distortion of Pascal graphs, approximating the Sierpinski gasket,
is bounded, which allows to obtain estimates for the convergence to zero of the
spectral gap as an application of the main result.
\end{abstract}


\end{frontmatter}

\section{Introduction}

Distortion of a finite metric space $X$ is, roughly speaking, a measure of how
well can $X$ be embedded into the Hilbert space or, more generally, a Banach space. 
The study of distortion of 
finite metric spaces has a long history, and in recent years the subject
received some attention, in part due to 
applications in theoretical computer science, where low distortion embeddings
imply good computational properties \cite{austin-naor-valette,
linial-london-rabinovich,linial-magen-naor}.

Another area of research which is on the border of discrete mathematics, probability,
algebra and computer science is the study of spectral, asymptotic 
and combinatorial properties of infinite families of graphs of constant, or uniformly
bounded, degree.  A characteristic  which is often at the heart of such studies is the 
first non-zero eigenvalue 
of the discrete Laplacian, known also as spectral gap.

Various relations between the spectral gap and other  natural characteristics of a  
finite graph are known, we refer to \cite{deabreu} for a
list of such relations. Two notable inequalities of this type are due 
to Alon and Milman \cite{alon-milman} and Chung \cite{chung}.
In this paper we study a similar relation, additionally involving the $\ell_p$-distortion of a finite graph.

Our first result  and the main tool in subsequent considerations is 
a general inequality relating the diameter, the first non-zero eigenvalue of the discrete 
$p$-Laplacian and  the $\ell_p$-distortion
of a family of graphs (Theorem \ref{theorem : main}). 
In order to prove it we introduce a certain mild regularity condition, called volume distribution. 
This condition guarantees that in some sense, sets of large 
volume cannot be small in diameter. It is of particular importance 
for infinite families of graphs, where we require uniform behavior of
the volume distribution. This condition is however natural and
it is satisfied by many families of graphs, examples of which are provided.

We apply the inequality to Cayley graphs of finite groups and
families of Schreier graphs with the set of vertices of the form $G/H_n$, where $G$ is a residually
finite group
generated by a finite set and $\set{H_n}_{n\in \NN}$ is a sequence of finite index
subgroups of $G$ with trivial intersection. 
The case of special interest is when the group $G$ is a finitely generated group of automorphisms 
of a regular rooted tree and $H_n=\operatorname{st}_G(u_n)$ is a 
stabilizer of a vertex $u_n$
of the $n$-th level, belonging to an infinite geodesic ray joining the root with infinity. The study of spectral properties of such graphs was initiated
in \cite{bartholdi-grigorchuk,grigorchuk-sunik-comptes,grigorchuk-sunik-standrews,grigorchuk-zuk}
and remains an active field of research.

One such family is obtained from the Hanoi Towers groups on 3 pegs.  These  
graphs, called Pascal graphs, were studied in \cite{grigorchuk-sunik-spectrum,quint} and 
geometrically  approximate the Sierpinski gasket. 
Our second result is  that the family of Pascal graphs has uniformly distributed volume 
(Theorem \ref{theorem : uniform volume distribution}) and
bounded $\ell_p$-distortion (Theorem \ref{theorem: distortion bounded}). 
As a consequence of this fact and our inequality we obtain
exponential convergence to zero for the spectral gap of the discrete $p$-Laplacian.
In the special case $p=2$ such convergence follows from \cite{grigorchuk-sunik-spectrum}.
In the process we also show that  information on distortion can give better estimates for 
convergence of  the spectral gap than some previously known inequalities.

It is worth mentioning that not all graphs having a fractal structure in the limit have
bounded distortion. Examples are furnished by the diamond graphs and Laakso graphs, which
are constructed using recursive substitution and whose $\ell_p$-distortion is unbounded
for any $p> 1$ \cite{johnson-schechtman}.

Finally, in the last section,
we discuss the remaining problems.

We would like to thank Zoran \v{S}uni\'{c} for carefully reading the paper and suggesting numerous improvements.

\section{Definitions} 

\subsection{Distortion}
Let $X$ be a finite metric space, $\mathcal{E}$ be a Banach space
 and let $f:X\to \mathcal{E}$ be a one-to-one map which is 1-Lipschitz:
$$\Vert f(x)-f(y)\Vert_\mathcal{E}\le d(x,y),$$
for all $x,y\in X$. Let
\begin{equation}\label{equation :  definition of L_f}
L_f=\max_{x,y\in X}\dfrac{d(x,y)}{\Vert f(x)-f(y)\Vert_\mathcal{E}}.
\end{equation}
\begin{definition}
The $\mathcal{E}$-distortion of $X$ is the number
\begin{equation}\label{equation : distortion defined by inf over L_f}
\disto_\mathcal{E}(X)=\inf L_f,
\end{equation}
where the infimum is taken over all $f:X\to \mathcal{E}$ which are one-to-one and 1-Lipschitz.
\end{definition}
Intuitively, distortion measures the most efficient way to embed $X$ into
$\mathcal{E}$.
By $\disto_{(p)}(X)$ we will denote the distortion of $X$ for $\mathcal{E}=\ell_p$.
In the case $p=2$ the number $\edisto(X)$ is 
called the \emph{Euclidean distortion} of $X$. Euclidean distortion is particularly 
well-studied, see for instance \cite{linial-london-rabinovich,linial-magen-naor,matousek} and the
references therein.

Bourgain \cite{bourgain} showed that Euclidean distortion of any finite metric space $X$
satisfies $c_{(2)}(X)\le C_B\log \vert X\vert$, where $\vert X\vert$ denotes the cardinality of $X$
and $C_B$ is a universal constant.
This bound is sharp and is realized by expander graphs (that is, infinite family of finite graphs
of uniformly bounded degree with spectral gap uniformly bounded away from zero) 
\cite{linial-london-rabinovich}.

\subsection{The discrete $p$-Laplacian and its eigenvalues} 
Let $\Gamma=(V,E)$ be a finite, undirected, connected graph with vertex set $V$ and edge 
set $E$. We write $x\sim y$ to denote the fact that vertices $x$ and $y$ are joined by an edge.
We allow loops 
and multiple edges and by $\omega(x,y)$ 
we denote the number of edges linking vertices $x$ and $y$.  
The set of vertices of $\Gamma$ can be viewed as a metric space when equipped 
with the combinatorial metric.

As usual,  given a set $S$, by $\ell_p(S)$ we denote the space of 
function $f:S\to \RR$ for which the $p$-norm 
$\Vert f\Vert_p=(\sum_{x\in S}\vert f(x)\vert^p)^{1/p}$ is finite. The symbol $1_A$
will denote the characteristic function of a set $A$.

Let $1< p<\infty$. The $p$-Laplacian $\Delta_p$ is an operator $\Delta_p:\ell_p(V)\to\ell_p(V)$, 
defined by the formula
$$\Delta_p f(x)=\sum_{x\sim y} (f(x)-f(y))^{[p]}\omega(x,y),$$
for $f:V\to \RR$,
where $a^{[p]}= \vert a\vert^{p-1} \sign(a)$.
For $p\neq 2$ the $p$-Laplacian is a non-linear operator, while for $p=2$
it is the standard linear discrete Laplacian.
The $p$-Laplacian is a well-known operator in the study of partial differential equations,
for graphs it was considered in e.g. 
\cite{amghibech,p-laplacian-1,buhler-hein,takeuchi}.

A real number $\lambda$ is an eigenvalue of the $p$-Laplacian $\Delta_p$
if there exists a function $f:V\to \RR$ such that 
$$\Delta_pf=\lambda f^{[p]}.$$
The eignevalues of the $p$-Laplacian are quite difficult to 
compute in the case 
$p\neq2$, due to non-linearity of $\Delta_p$.
Define
\begin{equation}\label{equation : lambda_1}
\lambda_1^{(p)}(\Gamma)=\inf\left\{\dfrac{\sum_{x\in V}\sum_{y\sim x}\vert f(x)-f(y)\vert^p \omega(x,y) }
{\inf_{\alpha\in \RR} \sum_{x\in V}\vert f(x)-\alpha\vert^p}\right\}
\end{equation}
with the infimum taken over all $f:V\to \RR$ such that $f$ is not constant.
It was proved in \cite{buhler-hein}, by means of the variational principle,
that $\lambda_1^{(p)}$ is  the smallest positive eigenvalue of the discrete
$p$-Laplacian $\Delta_p$ or the \emph{$p$-spectral gap}. For $p=2$ this agrees with the 
definition of the first positive eigenvalue as the minimizer of the Rayleigh quotient 
$\langle \Delta_2f,f\rangle/\langle f,f\rangle$
over all $f$ which are orthogonal to the constant functions on $V$,
since for $f:V\to \RR$ satisfying $\langle f,1_V\rangle=0$ we have
$$\Vert f\Vert_2^2=\inf_{\alpha\in \RR}\Vert f- \alpha 1_V\Vert_2^2.$$
In graph theory $\lambda_1^{(2)}$ is often referred to
as algebraic connectivity, see for instance \cite{deabreu}.

For $p=1$ the definition of $\lambda_1^{(1)}$ alone still makes sense and we will use state our results 
for $p\ge 1$, even though the interpretation of $\lambda_1^{(p)}$ as an eigenvalue of $\Delta_p$ is valid only for $p>1$.

\section{Distortion and the spectral gap} 
Given a finite metric space $X$ we denote by $\diam(X)$ its diameter and by $\vert X\vert $ its cardinality.
We start by defining a constant $\rho_{\epsilon}$ which describes certain geometric features of
$X$.
\begin{definition}
Let $X$ be a finite metric space. 
Given $0<\epsilon<1$ define the constant  $\rho_{\epsilon}(X)\in [0,1]$,
called the volume distribution,
by the relation
$$\rho_{\epsilon}(X)=\min\setd{\dfrac{\diam(A)}{\diam(X)}}{A\subseteq X \text{ such that }\vert A\vert \ge \epsilon \vert X\vert}.$$
\end{definition}
When $X=V$ is the vertex set of a finite graph $\Gamma$, we use the notation 
$\rho_{\epsilon}(\Gamma)=\rho_{\epsilon}(V)$. 

The following theorem gives the inequality between the $p$-spectral gap, diameter and 
the $\ell_p$-distortion.

\begin{theorem}\label{theorem : main}
Let $\Gamma$  be a finite, connected graph of degree bounded by $k$
and let $1\le p<\infty$. Then for every $0<\epsilon<1$,
$$
\lambda_1^{(p)}(\Gamma)\le C \left(\dfrac{\disto_{(p)}(\Gamma)}{ \diam(\Gamma)}\right)^{p},
$$
where $C=C(k,\epsilon,\rho_{\epsilon},p)=\dfrac{k}{1-\epsilon} \left(\dfrac{2}{\rho_{\epsilon} (\Gamma)} \right)^p$.
\end{theorem}

\begin{proof}[Proof of Theorem \ref{theorem : main}]
Let $\set{e_n}_{n\in\NN}$ be the standard basis vectors in $\ell_p(\NN)$ and let
$F:\Gamma\to \ell_p(\NN)$, where $F(x)=\sum_{n\in \NN}e_n F_n(x)$ and $F_n:V\to \RR$ 
are coordinate functions,
be a  one-to-one, $1$-Lipschitz map. 
Since 
$$\inf_{\alpha\in \RR}\left\vert \Vert F_n\Vert_p-\Vert \alpha 1_V\Vert_p\right\vert\le \inf_{\alpha\in \RR}\Vert F_n-\alpha 1_V\Vert_p\le \Vert F_n\Vert_p,$$
we  conclude that  for every $n\in \NN$ there exists $\alpha_n\in \RR$ such that
$$\inf_{\alpha\in \RR}\Vert F_n-\alpha 1_V\Vert_p^p=\Vert F_n-\alpha_n1_V\Vert_p^p,$$
and moreover,
\begin{equation}\label{equation: estimate on p-norm of alphas}
\vert \alpha_n\vert \le\dfrac{2\Vert F_n\Vert_p}{\Vert 1_V\Vert_p}=\dfrac{2\Vert F_n\Vert_p}{\vert V\vert^{1/p}}.
\end{equation}

We have $\sum_{n\in \NN}\vert \alpha_n\vert ^p<\infty$.
Indeed, inequality (\ref{equation: estimate on p-norm of alphas}) yields
$$\sum_{n\in \NN}\vert  \alpha_n\vert^p\le \dfrac{2^p}{\vert V\vert}\sum_{n\in \NN}\sum_{x\in V} \vert F_n(x)\vert^p= \dfrac{2^p}{\vert V\vert} \sum_{x\in V}\Vert F(x)\Vert^p_p<\infty.$$

By virtue of the above estimate, the vector $w=\sum_{n\in \NN}e_n\alpha_n$ is an element of $\ell_p$ and 
we can define a new embedding  $f:V\to \ell_p$ by  shifting by $w$: $f(x)=\sum_{n\in \NN}e_n f_n(x)$,
where
$$f_n(x)=F_n(x)-\alpha_n.$$  
Then for every $n\in\NN$ the inequality 
\begin{equation}\label{equation : inequality for f_n with lambda}
\lambda_1^{(p)}\sum_{x\in V}\vert f_n(x)\vert^p\le \sum_{x\in V}\sum_{y\sim x}\vert f_n(x)-f_n(y)
\vert^p\omega(x,y),
\end{equation}
holds, by the definition of $\lambda_1^{(p)}=\lambda_1^{(p)}$ and the choice of $\alpha_n$.
Then, applying (\ref{equation : inequality for f_n with lambda}) coordinate-wise, we have 
\begin{eqnarray*}
\sum_{x\in V}\Vert f(x)\Vert_p^p&=&\sum_{n\in \NN}\sum_{x\in V} \vert f_n(x)\vert^p\\
&\le&\dfrac{1}{\lambda_1^{(p)}}\sum_{n\in\NN}\sum_{x\in V}\sum_{y\sim x} \vert f_n(x)-f_n(y)\vert^p\omega(x,y)\\
&\le&\dfrac{1}{\lambda_1^{(p)}}\sum_{x\in V}\sum_{y\sim x}\Vert f(x)-f(y)\Vert_p^p\omega(x,y)\\
&\le&\dfrac{1}{\lambda_1^{(p)}} \sum_{x\in V}\sum_{y\sim x} \omega(x,y)\\
&\le&\dfrac{k\vert V\vert}{\lambda_1^{(p)}},
\end{eqnarray*}
since $F$ is 1-Lipschitz and, therefore, $\Vert f(x)-f(y)\Vert_p=\Vert F(x)-F(y)\Vert_p \le 1$ whenever $x\sim y$.
Thus at least $\epsilon\vert V\vert$ of $x\in V$ satisfy
$$\Vert f(x)\Vert_p\le \left(\dfrac{k}{(1-\epsilon) \lambda_{1}^{(p)}}\right)^{1/p}.$$ 
Consequently, we can find two points $x_0,y_0\in V$ such that 
$$d(x_0,y_0)\ge \rho_{\epsilon}(\Gamma)\diam(\Gamma)$$
and at the same time
$$\Vert f(x_0)-f(y_0)\Vert_p\le 2\left(\dfrac{k}{(1-\epsilon)\lambda_{1}^{(p)}}\right)^{1/p}.$$
For the constant $L_F$ we obtain
\begin{eqnarray*}
L_F&=&\max_{x,y\in V}\dfrac{d(x,y)}{\Vert F(x)-F(y)\Vert_p}\\
&\ge&\dfrac{d(x_0,y_0)}{\Vert f(x_0)-f(y_0)\Vert_p}\\
&\ge&\dfrac{1}{2}\rho_{\epsilon}(\Gamma)\diam(\Gamma)\cdot 
\left(\dfrac{k}{(1-\epsilon)\lambda_{1}^{(p)}}\right)^{-1/p}\\
&\ge&\dfrac{1}{2}\rho_{\epsilon}(\Gamma)\diam(\Gamma) \left(\lambda_1^{(p)}\right)^{1/p} \left(\dfrac{k}{1-\epsilon}\right)^{-1/p}.
\end{eqnarray*}
Since the right hand side of the inequality is independent of $F$ we pass to the infimum
over all $F:\Gamma\to \ell_p$ which are one-to-one and 1-Lipschitz and the assertion follows.
\end{proof}
 
\subsection{An Alon-Milman-type inequality for the $p$-spectral gap}
Alon and Milman proved in \cite{alon-milman} that for every graph  $\Gamma=(V,E)$ with 
degree bounded by $k$ the inequality
\begin{equation}\label{inequality : Alon-Milman}
\diam(\Gamma)\le 2\sqrt{\dfrac{2k}{\lambda_1^{(2)}(\Gamma)}}(\log_2\vert V \vert),
\end{equation}
holds.

On the other hand,  Bourgain proved in \cite{bourgain} that the Euclidean 
distortion of any finite metric space $X$ satisfies
$\edisto(X)\le C_B\log_2 \vert X\vert$ for some universal constant $C_B>0$. 
In fact his techniques show that $$\disto_{(p)}(X)\le C_B^{(p)}(\log_2 \vert X\vert),$$
for a universal constant $C_B^{(p)}$, depending only on $p$.   
This estimate
together with Theorem \ref{theorem : main} yields a general 
inequality between the diameter of a graph, number of vertices 
and $\lambda_1^{(p)}$. 
For every $0<\varepsilon <1$ we have 
\begin{equation}\label{equation: Alon-Milman for p-spectral gap}
\diam(\Gamma)\le \left(C(k,\epsilon,\rho_{\epsilon},p)
^{1/p}C_B^{(p)} \right)   \dfrac{\log_2 \vert V\vert}
{\left(\lambda_1^{(p)}(\Gamma)\right)^{1/p}}.
\end{equation}
For fixed $p\ge 1$, $\epsilon>0$ and a family of graphs of degree bounded by $k$, if $\rho_{\epsilon}(\Gamma)$ is uniformly 
bounded away from zero, then the constants $C(k,\epsilon,\rho_{\epsilon},p)
$ are uniformly
bounded above and in that case
Theorem \ref{theorem : main}, 
together with the result of Bourgain, allows to recover and generalize to any $p$, 
 the inequality (\ref{inequality : Alon-Milman}), up to a multiplicative constant depending on the family.

\section{Schreier graphs and uniform volume distribution}
In this section we consider infinite  families of Cayley graphs and  
Schreier graphs of groups of automorphisms of trees
(see \cite{bartholdi-grigorchuk,grigorchuk-nekrashevych-sushchanski,
grigorchuk-sunik-standrews,nekrashevych} 
for background on this topic).
We will demonstrate on examples that  for sequences of graphs Theorem \ref{theorem : main}
often gives  results asymptotically
close to optimal and in some cases gives better results than some previously known inequalites. 
Recall that
given a group $G$ with a finite generating set $A$ and a subgroup $P\le G$, the Schreier graph
$\Gamma=\Gamma(G,P,A)$ consists of the set of vertices being in bijection with left cosets $gP$ 
and the set of edges $E=\setd{(gP,agP)}{a\in A\cup A^{-1}}$. The Cayley graph of $G$ is a particular case
when $P=\set{1}$.  Every regular graph of even degree $2m$ can be represented as a Schreier graph of the
free group $\mathbb{F}_m$ with respect to a certain subgroup, see \cite[IV.A.15]{delaharpe} for a discussion.

We are interested in studying sequences of Schreier graphs of the form $\Gamma_n=(G,P_n,A)$, where 
 $\{P_n\} $ is a descending sequence of 
finite index subgroups with trivial intersection. In other words, we 
are considering 
a covering sequence $\{\Gamma_n\}$ of finite Schreier graphs 
(i.e., where $\Gamma_{n+1}$ covers $\Gamma_n$).
Groups with such chains are residually finite and act naturally  and level transitively on infinite, spherically 
homogeneous rooted trees, where $P_n$ serve as stabilizers of vertices of level $n$ of the rooted trees 
\cite{bartholdi-grigorchuk-sunik,grigorchuk-horizons,grigorchuk-nekrashevych-sushchanski,nekrashevych}.
We will consider here three examples: finite lamplighter groups, the group of intermediate growth 
constructed in \cite{grigorchuk-80,grigorchuk-84} and the  Hanoi Towers group $H^{(3)}$ introduced in
\cite{grigorchuk-sunik-comptes,grigorchuk-sunik-standrews}. Many more example of such graphs
can be found in \cite{bartholdi-grigorchuk,grigorchuk-savchuk-sunik,kaimanovich,teplayev}.

We are interested in three characteristics of such sequences: diameters, $\lambda_1$ and distortion. 
It is easy to see that the asymptotic behavior of diameters and distortion does not depend on 
the generating
set $A$, however it is not clear if the same is true for $\lambda_1$, therefore we will always indicate the
generating set we are using.

\subsection{Uniform volume distribution}
We now want to consider families of graphs for which  the constant 
$C(k,\epsilon,\rho_{\epsilon},p)
 $ is uniformly bounded above. This amounts to 
controlling $\rho_{\epsilon}$ and motivates the following
notion.
 \begin{definition}\label{definition : uniform volume distribution}
 Given a family of finite
metric spaces $\mathcal{X}=\set{X_i}$ and $0<\epsilon<1$ we say that the family 
$\mathcal{X}$ has uniform volume distribution if there exist $\epsilon>0$ and $c>0$ such that
$\rho_{\epsilon}(X_i)\ge c$
for every $X_i\in \mathcal{X}$.
\end{definition}

Such families will be of particular interest.
Of course not all families of bounded degree graphs have uniform volume distribution.
\begin{example}\normalfont
Consider a $3$-regular tree $T$ and let $B(n)$ denote a ball of radius $n$ around a 
fixed point $x_0$. 

Define the graph $\Gamma_n$ as $B(n)$ with a path $P_n$ of length $\vert B(n)\vert$ attached at $x_0$.  
Then $\vert B(n)\vert =\vert \Gamma_n\vert/2$, however $\delta(B(n))=2n$ and
$\delta(\Gamma_n)\ge \vert P_n\vert\ge \vert B(n)\vert$, where the latter grows exponentially, 
so that the ratio 
$\delta(B(n))/\delta(\Gamma_n)$ tends to 0.
\end{example}
Fortunately, many families satisfy the conditions of Definition 
\ref{definition : uniform volume distribution} in a natural way. 
Recall that a graph is vertex transitive if the group of automorphisms acts transitively on its vertices.
Cayley graphs of groups are examples of vertex transitive graphs.
\begin{proposition}
Let $\Gamma$ be a vertex transitive connected graph on at least 3 vertices. 
Then $\rho_{1/2}(\Gamma)\ge \dfrac{1}{4}$.
\end{proposition}
\begin{proof}

By $B_x(r)=\setd{y}{d(x,y)\le r)}$ we denote the ball or radius $r$ centered at $x$.
If $\diam(\Gamma)\le 2$ then the claim can verified directly.
Assume that $\diam(\Gamma)\ge 3$.
First let us prove  that $\left\vert B_x(\diam(\Gamma)/4)\right\vert <\vert V\vert/2$ for any vertex $x\in V$.  Assume the contrary.
Take a pair of points $x,y\in V$ such that $$d(x,y)=\diam(\Gamma).$$
We have
$$ B_x(\diam(\Gamma)/4)\cap B_y(\diam(\Gamma)/4)=\emptyset$$
and there is a point which does not belong to either of the two balls.
Indeed, if there is no such point then 
$$\delta(\Gamma)=d(x,y)\le 2(\delta(\Gamma)/4)+1=\delta(\Gamma)/2+1,$$
which implies $\delta(\Gamma)\le 2$. 
By the triangle inequality and homogeneity we obtain  a contradiction
$$\vert V\vert> 2\cdot \vert B_x(\diam(\Gamma)/4)\vert \ge 2\cdot\frac{\vert V\vert}{2}\ge\vert V\vert.$$

Now assume that $\rho_{1/2}(\Gamma)<\dfrac{1}{4}$. Then there exists a set $A\subseteq V$ such that 
$$\vert A\vert \ge \dfrac{\vert V\vert}{2} \ \ \ \ \text{and}\ \ \ \ \ \  \diam(A)<\dfrac{\diam(\Gamma)}{4}.$$
Choose any point $x_0\in A$. Then $A\subseteq B_{x_0}(\diam(\Gamma)/4)$.
Thus $$\vert A\vert \le \left\vert B_{x_0}(\diam(\Gamma)/4)\right\vert<\dfrac{\vert V\vert}{2}$$
and we again get a contradiction.
\end{proof}

\begin{corollary}
The family of all vertex transitive connected graphs on at least three vertices 
has $1/4$-uniform volume distribution.
\end{corollary}

We recall an inequality  due to Chung \cite{chung}, who showed that for a $k$-regular graph 
\begin{equation}\label{inequality : Chung}
\diam(\Gamma)\le \left\lceil \dfrac{\log(\vert V \vert-1)}{\log k/\alpha}\right\rceil,
\end{equation}
where $\alpha=\vert \alpha_2\vert$ for the eigenvalues of the adjacency matrix
$\alpha_1,\alpha_2\dots$ satisfying $\vert \alpha_1\vert\ge \vert \alpha_2\vert\ge\dots\ge \alpha_n$.
Below we will compare some of our results with those which follow from inequality 
(\ref{inequality : Chung})
Note that in comparison with inequality (\ref{inequality : Alon-Milman}), Theorem \ref{theorem : main} automatically gives better convergence as soon as a family of graphs has
distortion better than $O(\log\vert V_i\vert)$. Indeed, inequality (\ref{inequality : Alon-Milman})
follows from Theorem \ref{theorem : main} and Bourgain's upper bound $C_B^{(2)}\log(\vert V\vert)$ on
the distortion.
A version of  inequality (\ref{inequality : Chung}) in the case of directed graphs, including
Cayley graphs of finite groups, was also studied recently in 
\cite{chung-directed}. We direct the reader to \cite{chung-book} for details.
\subsection{Finite wreath products}
Consider the wreath product $\ZZ_2\wr \ZZ_n$  of cyclic groups
of order $2$ and $n$, that is the semidirect product
$(\oplus_{i=1}^n \ZZ_2)\rtimes \ZZ_n$, where the action of $\ZZ_n$ on $\oplus_{i=1}^n\ZZ_2$
is given by a coordinate shift. 
The natural generating set of the wreath product $\ZZ_2\wr\ZZ_n$ is 
$$\set{(1,0,0,\dots,0),a},$$
where $a$ is the generator of $\ZZ_n$, and
one can easily prove that the diameters of the Cayley graphs with respect to this system of generators
grow linearly. Spectra of finite wreath products as above were studied in 
\cite{grigorchuk-zuk-dedicata}. It was shown in \cite{austin-naor-valette}
that the Euclidean distortion in this case is $O(\sqrt{\log n})$. 
Therefore, by Theorem \ref{theorem : main} and uniform volume distribution, 
we have
\begin{proposition} There exists a constant $C>0$ such that
\begin{equation}\label{equation : estimate for wreath products}
\lambda_1^{(2)}(\ZZ_2\wr\ZZ_n)\le C\left(\dfrac{\sqrt{\log n}}{n}\right)^2=C\dfrac{\log n}{n^2}.
\end{equation}
\end{proposition}
In comparison,  inequality (\ref{inequality : Chung}) gives 
$$\lambda_1^{(2)}(\ZZ_2\wr\ZZ_n)\le 1-(n^2-1)^{(-1/n)}.$$
Then one can verify that the estimate (\ref{equation : estimate for wreath products}) 
is asymptotically stronger.

\subsection{The group of intermediate growth}
\begin{figure}\includegraphics[width=250pt]{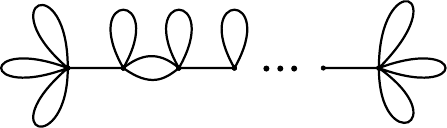}
\caption{The Schreier graph at level $n$ of the group of intermediate growth has $2^n$ vertices}\label{figure : Schreier graph of the group of intermediate growth}
\end{figure}

Consider the group $G$ of intermediate growth generated by automorphisms $a,b,c,d$ of order 2 
of a binary rooted tree. The description of this group can be found in many places, for instance
\cite{grigorchuk-80,grigorchuk-84,grigorchuk-pak}, while the spectral properties of Schreier 
graphs of $G$ were studied in \cite{bartholdi-grigorchuk}.
Without getting into details let us mention that the shape of the Schreier graphs 
is shown in Figure \ref{figure : Schreier graph of the group of intermediate growth}.

We have in this case
$\diam(\Gamma_n)=2^{n}-1$. 
It is easy to see that the $\ell_p$-distortion of $\Gamma_n$ is bounded and that $\Gamma_n$ have
uniform volume distribution with $\rho_{1/2}(\Gamma_n)\ge 1/2$.
Thus we obtain from Theorem \ref{theorem : main},
\begin{proposition}
The $p$-spectral gap of $\Gamma_n$ converges exponentially to 0.
More precisely, there exists a constant $C>0$ such that
$$\lambda_1^{(p)}(\Gamma_n)\le \dfrac{C}{2^{np}}.$$
\end{proposition}
In the case $p=2$ the spectrum of $\Gamma_n$ was precisely computed in \cite{bartholdi-grigorchuk} 
and gives the asymptotics
$C_14^{-n}\le \lambda_1(\Gamma_n)\le C_24^{-n}$.
In comparison, inequality (\ref{inequality : Chung}) gives
$\lambda_1^{(2)}(\Gamma_n)\le C_3 n2^{-n}$. Here 
$C_1,C_2,C_3$ are positive constants.

\section{Hanoi Towers group on 3 pegs and its Schreier graphs}
The Hanoi Towers game on $k\ge 3$ pegs leads to a famous combinatorial problem 
\cite{hinz,szegedy},
which remains unsolved for $k\ge 4$ (although it is solved asymptotically in \cite{szegedy}).
The case $k=3$, despite its easy solution, attracted the attention of many researchers
because of its relations to different topics in mathematics, including  studies around
the Sierpi\'{n}ski gasket \cite{grigorchuk-sunik-spectrum,quint,teplayev}.

In \cite{grigorchuk-sunik-comptes} it was discovered that the problem has an interesting 
connection to algebra,
namely to the theory of self-similar groups, branch groups and iterated monodromy
groups.

\begin{figure}[!ht]
\begin{center}
\includegraphics[width=200pt]{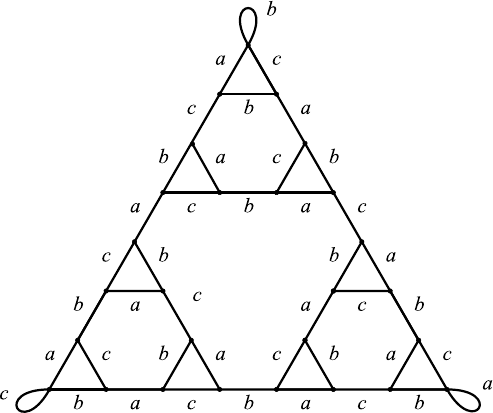}
\caption{The Schreier graph of $H^{(3)}$ at level 3}\label{sierpinski3}
\end{center}
\end{figure}

The natural sequence of graphs related to the game on 3 pegs is a sequence of the so-called
Pascal graphs, which is closely related to the sequence of Sierpi\'{n}ski graphs studied in fractal theory 
\cite{teplayev}. Slight modifications of these graphs also arise as Schreier graphs for the natural 
action of the 
associated Hanoi Tower group $H^{(3)}$ on the ternary tree. Their shape is shown in 
Figure \ref{sierpinski3}. Let us briefly define the corresponding group
and the sequence of graphs. 
We direct the reader to 
\cite{grigorchuk-nekrashevych-sushchanski,grigorchuk-sunik-comptes,grigorchuk-sunik-standrews,
grigorchuk-sunik-spectrum} for background on the Hanoi Towers groups. 

Let $\mathcal{T}_3$ denote the 3-regular rooted tree. The set of vertices 
of $\mathcal{T}_3$ can be viewed as 
the set of elements of the free monoid $\mathcal{A}^*$ of finite words over the alphabet 
$\mathcal{A}=\set{0,1,2}$. 
Given a permutation $\sigma\in \Sigma_3$  in the symmetric group define the automorphism 
$a_{\sigma}$ of $\mathcal{T}_3$ by
the recursive rule
$$a=\sigma(a_0,a_1,a_2),$$
where $a_i$ is the identity automorphism if $i$ is in the support of $\sigma$ and $a_i=a$
otherwise. The automorphisms $a_{(ij)}$ act on the set of vertices $\mathcal{T}_3$ by the following 
recursive formulas:
\begin{eqnarray*}
a_{(ij)}(iw)&=&jw\\
a_{(ij)}(jw)&=&iw\\ 
a_{(ij)}(xw)&=&xa_{(ij)}(w)
\end{eqnarray*}
for $x\neq {i,j}$, where $i,j\in \mathcal{A}$ and $w\in \mathcal{A}^*$. 
Denote $a=a_{(01)}$, $b=a_{(02)}$ and $c=a_{(12)}$.
The Hanoi Towers group on $3$ pegs is the group 
$H^{3}\subseteq \operatorname{Aut}(\mathcal{T}_3)$ generated by the
automorphisms $a$, $b$ and $c$.
Consider now the subgroup $P_n=\operatorname{st}_{H^{(3)}}(1^n)\le H^{(3)}$. 
The Schreier graph of this subgroup is shown in Figure \ref{figure : H_n+1,3}.
The spectral theory of Pascal graphs was studied in 
\cite{grigorchuk-sunik-standrews} and \cite{quint}, from which it follows that
 $\lambda_1\le C 5^{-n}$. We will obtain information about the distortion 
 of Schreier graphs (which are graphs of actions on levels) and apply Theorem \ref{theorem : main}.

\begin{theorem}\label{theorem : uniform volume distribution}
The graphs $\Gamma_n$ have uniform volume distribution. More precisely, 
$\rho_{2/3}({\Gamma_n})\ge 1/2$ for every $n\in \NN$.
\end{theorem}
\begin{proof} 
We have $\diam(\Gamma_n)=2^{n}-1$.
Denote by $S_1$, $S_2$ and $S_3$ the three copies of 
${\Gamma_{n-1}}$ which
are naturally embedded in ${\Gamma}_{n}=(V,E)$, as in Figure \ref{figure : H_n+1,3}.
Consider a set $A\subseteq V$ 
which satisfies $\vert A\vert \ge 2/3 \vert V\vert$. Then 
we consider two cases.

If $\vert A\vert=2/3\vert V\vert$ and $A=S_i\cup S_j$ for some $i,j$ 
then  $\diam(A)=\diam (\Gamma_{n})$.

Let $A$ intersect all three $S_i$, $i=1,2,3$.
and let $x_i\in A$, $i=1,2,3$ belong to $S_i$. Consider two cases. 
The first one is when the distance between two of the points, say $x_2$ and $x_3$ is realized by
a geodesic that passes through $S_1$. In that case the distance $d(x_3,x_2)\ge \delta(\Gamma_{n-1})+2\ge \delta(\Gamma_n)/2$. 

The second case is when the distance between $x_i$ and $x_j$ is not realized by a geodesic through $S_k$
where $i,j,k\in \{1,2,3\}$ and are all different.
For each $i=1,2,3$ let $v_i^j$ for  $j=1,2$ denote
the two vertices of  the triangle $S_i$ that are not vertices of the large triangle, as shown in Figure \ref{figure : H_n+1,3}. 
Then we have
\begin{eqnarray*}
\diam(A)&\ge&\max\set{d(x_1,x_2), d(x_2,x_3),d(x_3,x_1)}\\
&\ge& \dfrac{1}{3}\left(d(x_1,x_2)+d(x_2,x_3)+d(x_3,x_1)\right)\\
&\ge& \dfrac{1}{3} \left(d(x_1,v_1^2)+d(v_2^1,x_2)+1\right. \\
&&\ \ \ \ \ \ \ \ \left.+d(x_2,v_2^2)+d(v_3^1,x_3)+1+d(x_3,v_3^2)+d(v_1^1,x_1)+1\right)\\
&\ge&\dfrac{1}{3}\left(3+ d(v_1^1,v_1^2)+d(v_2^1,v_2^2)+d(v_3^2,v_3^3)\right)\\
&\ge&1+\diam(\Gamma_{n-1})\ge\dfrac{\diam(\Gamma_n)}{2},
\end{eqnarray*}
which finishes the proof.
\begin{figure}
\includegraphics[width=200pt]{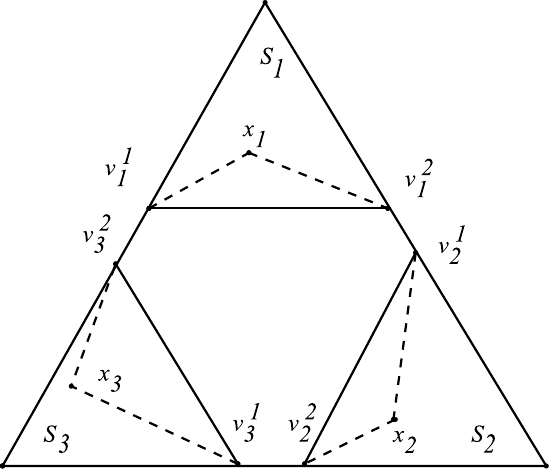}
\caption{The graph ${\Gamma}_{n}$}\label{figure : H_n+1,3}
\end{figure}
\end{proof}

\begin{theorem}\label{theorem: distortion bounded}
For any $p\ge 1$ the $\ell_p$-distortion of the sequence $\{\Gamma_n\}$ is 
uniformly bounded above.
\end{theorem}
Before we proceed with the proof we recall the definition of a quasi-isometry,
which is a standard notion of metric equivalence in 
coarse geometry and geometric group theory.
\begin{definition}
Let $X$ and $Y$ be metric spaces. A map $f:X\to Y$ is a quasi-isometry if there
exist constants $L\ge 1$ and $K,C\ge 0$ such that 
$$\dfrac{1}{L}d_X(x,y)-C\le d_Y(f(x),f(y))\le Ld_X(x,y)+C$$
and the image of $X$ is $K$-dense in $Y$; that is, for every $y\in Y$ there exists 
$x\in X$ such that $d_Y(f(x),y)\le K$.
\end{definition}

\begin{figure}
\includegraphics[height=200pt]{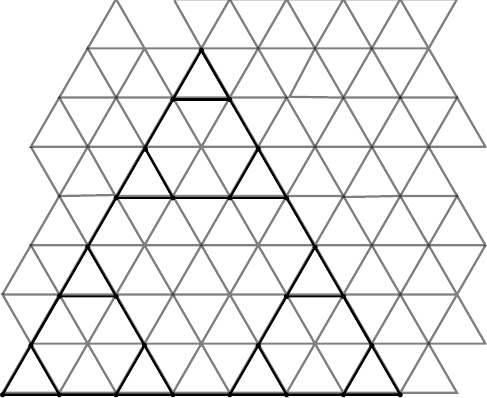}
\caption{Embedding $\Gamma_n$ into $\RR^2$.}\label{figure : triangular net on the plane}
\end{figure}

\begin{proof}[Proof of Theorem \ref{theorem: distortion bounded}]
We first consider the case $p=2$. Take the triangulation of the plane by equilateral triangles,
as in Fig. \ref{figure : triangular net on the plane}. The set of vertices of this net will
be called $M$.
This set can be equipped with the combinatorial metric on the 1-skeleton $N$
of the triangulation. With this metric, $N$ is quasi-isometric to the plane equipped 
with any $p$-norm, $p\ge 1$.
Indeed, $M$ can be viewed as the Cayley graph of the group $\ZZ^2$ with the generating
set $$\set{(1,0),(-1,0),(0,1),(0,-1),(-1,1),(1,-1)}.$$
Theorem  \ref{theorem: distortion bounded} will be proved once we show that the graphs ${\Gamma_n}$  
are subgraphs of $(M,N)$ with the sets of vertices embedded
quasi-isometrically
with uniform constants into $M$.

First, note that ${\Gamma}_1$ embeds into $M$ with distortion 
less than $3$.
We will prove by induction on $n$ that ${\Gamma}_{n}$ embeds into $M$ with distortion bounded by $3$.

 Denote by $S_1,S_2,S_3$ the three copies of 
${\Gamma}_{n-1}$ embedded in ${\Gamma_n}$, starting with $S_1$ being the top one and continuing
in the clockwise direction, as in Figure \ref{figure : H_n+1,3}. 
If $x,y$ are verticies in $S_i$ for some $i=1,2,3$ then the claim follows from the inductive assumption.
It suffices to prove the estimate for $x\in S_1$ and $y\in S_2$, the other cases being completely
analogous. We consider the infinite straight line $l$ passing through $x$, parallel to that side 
 ${\Gamma_n}$ (viewed as a triangle on the plane)
which intersects both $S_1$ and $S_2$. The point $y\in S_2$ can be on the left side of 
$l$ or on the right side of $l$ (see Figures \ref{figure : y to the left of x} and \ref{figure : y to the right of x}). We consider each case separately.\\
\begin{figure}
\includegraphics[width=200pt]{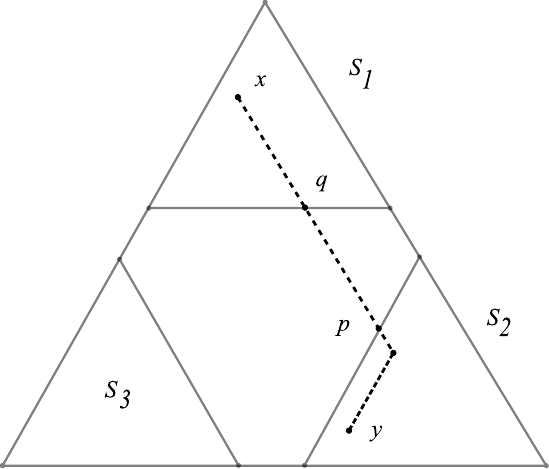}
\caption{$y$ is to the left of $x$}\label{figure : y to the left of x}
\end{figure}

\noindent\emph{Case 1:  $y$ is to the left of $l$.} Consider the geodesic linking $x$ and $y$ in $M$ which consists
of two straight segments: the first segment is contained in $l$; the other segment is the 
segment of a straight line through $y$ parallel to the side of ${\Gamma_n}$ which does not intersect 
$S_2$. Let $p$ be the point of intersection of this geodesic with that side of $S_2$
which is the only side not contained in any side of ${\Gamma}_{n}$.
Then the geodesic segment $[y,p]$ is completely contained in $S_2$ and the estimate for
distortion for $\Gamma_{n-1}$ applies.
Let $q$ denote the point of intersection of the geodesic with that side of $S_1$ which is the only side not contained in any side of ${\Gamma}_{n}$. Then the geodesic segment $[x,q]$ lies 
completely inside of $S_1$ and again the estimate for  distortion $\Gamma_{n-1}$ applies. 

For the part of the geodesic linking $p$ and $q$ we note that it lies on the line parallel
to the side of ${\Gamma_n}$ and so we have $$d_M(p,q)\ge \dfrac{1}{3}d_{\Gamma_n}(p,q).$$
This estimate follows from the fact that the geodesic linking $p$ and $q$ in $M$ can be viewed 
as the side of a equilateral triangle of side length $d_M(p,q)$. Then
$d_{\Gamma_n}(p,q)\le 2d_M(p,q)+1\le 3d_M(p,q)$.
Summarizing, we have
\begin{eqnarray*}
d_M(x,y)&=&d_M(x,p)+d_M(p,q)+d_M(q,y)\\
&\ge&\dfrac{1}{3}d_{\Gamma_n}(x,q)+\dfrac{1}{3}d_{\Gamma_n}(q,p)+\dfrac{1}{3}d_{\Gamma_n}(p,y)\\
&\ge&\dfrac{1}{3}d_{\Gamma_n}(x,y).
\end{eqnarray*}

\noindent\emph{Case 2: $y$ is to the right of $l$.}
Consider the geodesic $\gamma$ linking $x$ and $y$, consisting of a segment of $l$ and
a segment of a horizontal straight line between $y$ and the intersection with $l$. 
Let $v$ denote the point of intersection of these two lines.
This leads to two subcases:\\
\begin{figure}
\includegraphics[width=200pt]{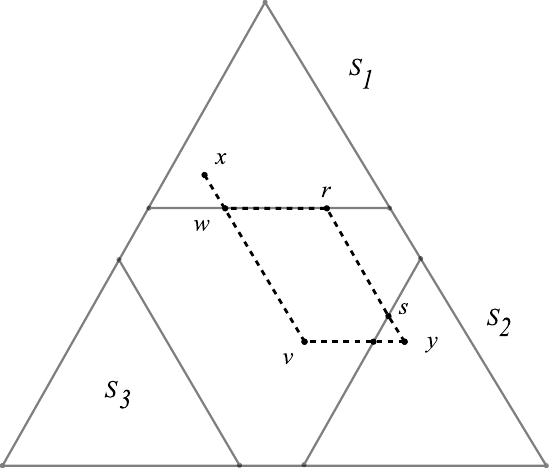}
\caption{$y$ is to the right of $x$}\label{figure : y to the right of x}
\end{figure}

\noindent a) If $v$ lies inside of $S_2$ we use arguments as above.\\

\noindent b) If $v$ lies outside of 
$S_2$ we replace the geodesic $\gamma$ by another geodesic as follows. 
We take the point of intersection $w$ of 
the geodesic $\gamma$ with the boundary of $S_1$. 
The points $w$ and $y$ can now be linked by geodesic
$\gamma'$ which consists of two straight line segments:
a segment of the side of $S_1$ containing $w$ and a segment
of a straight line through $y$ parallel to the side of $\Gamma_n$ intersecting 
both $S_1$ and $S_2$.
Our new geodesic 
is the segment $d[x,w]\cup \gamma' $, see Figure \ref{figure : y to the right of x}.

Consider now the new geodesic linking $x$ and $y$ to be oriented in the direction from $x$ to $y$.
Let $r$ denote the last point of intersection of the geodesic with $S_1$ and 
let $s$ denote the first point of intersection of the geodesic with $S_2$. 
 Then again, the segment linking $r$ and $s$ is on a line parallel to one side
 of ${\Gamma_n}$ while the segments from $x$ to $r$ and from $s$ to $y$ are completely
 contained in $S_1$ and $S_2$ respectively. 
 Thus we have an estimate:

 \begin{eqnarray*}
d_M(x,y)&=&d_M(x,r)+d_M(r,s)+d_M(s,y)\\
&\ge&\dfrac{1}{3}d_{\Gamma_n}(x,r)+\dfrac{1}{3}d_{\Gamma_n}(r,s)+\dfrac{1}{3}d_{\Gamma_n}(s,y)\\
&\ge&\dfrac{1}{3}d_{\Gamma_n}(x,y).
\end{eqnarray*}
This proves that the $\ell_p$-distortion of $\Gamma_n$ is bounded.
\end{proof}

Applying the above facts together with our inequality we obtain
\begin{theorem} 
The $p$-spectral gap of the graphs $\Gamma_n$ converges exponentially to zero for $1\le p<\infty$.
More precisely, there exist a constant $C>0$ (depending on $p$) such that for every $n\in \NN$ we have
$$\lambda_1^{(p)}({\Gamma_n}) \le C2^{-np}.$$
\end{theorem}
Thus for $p=2$ we have $\lambda_1\le C4^{-n}$.
The exponential convergence in the case $p= 2$ follows from 
\cite{grigorchuk-sunik-spectrum}

\begin{remark}\normalfont
The Sierpinski graphs \cite{teplayev}  are very similar to the graphs $\Gamma_n$ above, 
the difference being the three copies of $\Gamma_{n-1}$ in $\Gamma_n$ are linked at  the
vertices, without the connecting edges. The techniques used in this section give conclusions
 similar as above for the Sierpinski graphs.
\end{remark}

\section{Remaining questions}
There are some problems that we consider interesting in this context.
The Basilica group was introduced in \cite{grigorchuk-zuk} and its Schreier graphs were studied
in \cite{bondarenko,bartholdi-grigorchuk-nekrashevych,d'angeli}. These graphs can be constructed by a replacement
algorithm and have a tree-like structure with cycles of increasing diameters.

\begin{problem}
What is the distortion of Schreier graphs of the Basilica group? 
\end{problem}
We conjecture that the distortion of these Schreier graphs is unbounded.
The spectra of Laplacians on the (slightly modified) Schreier graphs of the Basilica group
were computed in \cite{rogers-teplyaev}, while for Schreier graphs associated to Hanoi Tower
groups $H^{(k)}$, $k\ge 4$, the spectra are not known. In the latter case it is known that the diameters
$\diam(\Gamma_n)$ grow asymptotically as $\exp(n^{1/(k-2)})$ and Chung's inequality gives 
$$\lambda_1^{(2)}(\Gamma_n)\le Cn\exp(-n^{1/(k-2)}).$$

The Hanoi Tower group we considered corresponds to the case $k=3$ pegs. 
The question of computing the spectrum of Schreier graphs associated to Hanoi Tower groups
on $k\ge 4$ pegs was posed in \cite{grigorchuk-sunik-standrews}. Here we state
\begin{problem}
What is the asymptotic behavior of $\lambda_1^{(2)}(\Gamma_n)$ for the Schreier graphs associated
 to the Hanoi
Tower groups $H^{(k)}$ for $k\ge 4$? 
\end{problem}

The possibility of applying Theorem \ref{theorem : main} to the above problem also yields 

\begin{problem}
What is the distortion of the family of Schreier graphs $\{H_n^{(k)}\}_{n\in \NN}$ when $k\ge 4$?
\end{problem}

\end{document}